\documentclass[a4paper, 12pt]{amsart}

\usepackage{amsfonts}
\usepackage{amsmath}

\newtheorem{theorem}{Theorem}[section]
\newtheorem{lemma}[theorem]{Lemma}
\newtheorem{proposition}[theorem]{Proposition}
\newtheorem{corollary}[theorem]{Corollary}
\theoremstyle{definition}
\newtheorem{definition}[theorem]{Definition}

\begin{document}
\date{2010-5-15}
\title[Joinings and subsystems]{Relatively independent joinings and
subsystems of W*-dynamical systems}
\author{Rocco Duvenhage}
\address{Department of Physics\\
University of Pretoria\\
Pretoria 0002\\
South Africa}
\email{rocco.duvenhage@up.ac.za}

\begin{abstract}
Relatively independent joinings of W*-dynamical systems are constructed.
This is intimately related to subsystems of W*-dynamical systems, and
therefore we also study general properties of subsystems, in particular
fixed point subsystems and compact subsystems. This allows us to obtain
characterizations of weak mixing and relative ergodicity, as well as of
certain compact subsystems, in terms of joinings.
\end{abstract}

\maketitle

\section{Introduction}

In the study of the general mathematical structure of quantum dynamical
systems and quantum statistical mechanics, the operator algebraic approach
has proven very valuable. In particular the framework of von Neumann
algebras, in which case the dynamical system is called a W*-dynamical system
(see Section 2 for a precise definition), provides a natural arena for
ergodic theory which extends the classical measure theoretic framework.
Refer to \cite{BR, BR2} for an account of many aspects of these topics.

The first aim of this paper is to study the extension of the concept of a
relatively independent joining of two dynamical systems in classical ergodic
theory to the noncommutative framework of W*-dynamical systems. It is
essentially a generalized way of forming the tensor product of two systems
which takes into account a common subsystem of the systems. A clear
exposition of the classical case can be found in \cite[Chapter 6]{G}.

Since the early works \cite{F67, Rud}, joinings have been a useful tool in
classical ergodic theory (again see \cite{G}). This paper, building on \cite{D, D2},
is part of a programme to systematically develop the theory of
joinings in the noncommutative case. The goal is to eventually have a
similarly useful tool for noncommutative dynamical systems.

Joinings have in fact already gradually found some use in noncommutative
dynamical systems, in particular related to dynamical entropy \cite{ST}
where a special case of joinings appears (also see \cite[Section 5.1]{NS}).
Early work related to joinings in a noncommutative setting goes back to at
least \cite{A}, where disjointness was studied in the context of a
noncommutative extension of topological dynamics (see in particular
\cite[Definition 2.9 and Section 5]{A}).

Since subsystems form an integral part of relatively independent joinings,
we also study the properties of particular subsystems, namely the fixed
point subsystem and the compact subsystem generated by the eigenoperators of
the dynamics. These subsystems are related to the properties of W*-dynamical
systems, namely to ergodicity and weak mixing respectively. For subsystems
more generally we will see in Section 2 and onward that the modular group of
the state of the W*-dynamical system in question plays an essential role in
the definition and application of subsystems relevant to relatively
independent joinings. This fits in very naturally with the physically
relevant case in which the dynamics itself is given by the modular group of
the state, but we will see that the framework is much more widely
applicable. Note that subsystems have already proven useful in the study of
noncommutative dynamical systems; see for example \cite[Sections 3 and 4]{AET}
and \cite[Definition 2.9 and Lemma 2.10]{P} for recent work.

After summarizing the basic framework in Section 2, we construct relatively
independent joinings of two W*-dynamical systems in Section 3 and derive a
few useful facts regarding them. In Sections 4 and 5 we then study relative
ergodicity and compact subsystems respectively, to illustrate how relatively
independent joinings fit into the theory of W*-dynamical systems. In the
latter two sections the fixed point subsystem and compact subsystem
generated by the eigenoperators respectively play a central role. In Section
6 we conclude the paper by studying further properties of this compact
subsystem when the W*-dynamical system is ergodic.

\section{Basic definitions, notations and background}

We use the same basic definitions as in \cite{D, D2}. For convenience we
summarize them here, along with some additional definitions. Simultaneously
this fixes notations that will be used throughout the rest of the paper.
Some related background material, in particular regarding Tomita-Takesaki
theory (or modular theory) is also discussed. A general notation that we use
often is $B(X)$ to denote the space of all bounded linear operators
$X\rightarrow X$ on a normed space $X$. The identity element of a group will
be indicated by $1$. In the remainder of this paper W*-dynamical systems are
referred to simply as ``systems'' and they are defined as follows:

\begin{definition}
A \emph{system} $\mathbf{A}=\left( A,\mu ,\alpha \right) $ consists
of a faithful normal state $\mu $ on a (necessarily $\sigma$-finite)
von Neumann algebra $A$, and a representation $\alpha :G\rightarrow
$ Aut$(A):g\mapsto \alpha _{g}$ of an arbitrary group $G$ as
$\ast$-automorphisms of $A$, such that $\mu \circ \alpha _{g}=\mu $
for all $g$. We call the system $\mathbf{A}$ an \emph{identity
system} if $\alpha _{g}=$ id$_{A}$ for all $g$ where
id$_{A}:A\rightarrow A$ is the identity mapping, while we call it
\emph{trivial} if $A=\mathbb{C}1_{A}$ where $1_{A}$ (often denoted
simply as $1$) is the unit of $A$.
\end{definition}

In the rest of the paper the symbols $\mathbf{A}$, $\mathbf{B}$, $\mathbf{F}$
and $\mathbf{R}$ will denote systems $\left( A,\mu,\alpha \right)$,
$\left(B,\nu ,\beta \right)$, $\left(F,\lambda ,\varphi \right) $ and
$\left(R,\psi,\rho\right)$ respectively,\ all making use of actions of the
same group $G$.

\begin{definition}
A \textit{joining} of $\mathbf{A}$ and $\mathbf{B}$ is a state $\omega$
(i.e. a positive linear functional with $\omega(1)=1$) on the algebraic
tensor product $A\odot B$ such that $\omega\left(a\otimes1_{B}\right)
=\mu(a) $, $\omega\left(1_{A}\otimes b\right) =\nu(b)$ and
$\omega\circ\left(\alpha_{g}\odot\beta_{g}\right)=\omega$ for all $a\in A$,
$b\in B$ and $g\in G$. The set of all joinings of $\mathbf{A}$ and $\mathbf{B}
$ is denoted by $J\left(\mathbf{A},\mathbf{B}\right)$. We call $\mathbf{A}$
\textit{disjoint} from $\mathbf{B}$ when $J\left(\mathbf{A},\mathbf{B}
\right) =\left\{\mu\odot\nu\right\}$.
\end{definition}

We will also have occasion to use a more general concept, namely if $A$ and
$B$ are von Neumann algebras with faithful normal states $\mu $ and $\nu $
respectively, then a \emph{coupling} of $\left(A,\mu \right)$ and $\left(
B,\nu \right)$ is a state $\omega$ on $A\odot B$ such that $\omega \left(
a\otimes 1_{B}\right) =\mu (a)$ and $\omega \left( 1_{A}\otimes b\right)
=\nu (b)$.

The modular group of a faithful normal state $\mu $ on a von Neumann algebra
$A$ will be denoted by $\sigma ^{\mu }$, and its elements by $\sigma
_{t}^{\mu }$ for every $t\in \mathbb{R}$. The cyclic representation of $A$
obtained from $\mu $ by the GNS construction will be denoted by $\left(
H_{\mu },\pi _{\mu },\Omega _{\mu }\right) $; this notation will in fact
also be used in the case of an arbitrary state on a unital $\ast $-algebra.
The associated modular conjugation will be denoted by $J_{\mu }$, and we let
\begin{equation*}
j_{\mu }:B(H_{\mu })\rightarrow B(H_{\mu }):a\mapsto J_{\mu}a^{\ast}J_{\mu}
\end{equation*}
for which we note that $j_{\mu }^{-1}=j_{\mu}$. We will also use the
notation
\begin{equation*}
\gamma _{\mu }:A\rightarrow H_{\mu }:a\mapsto \pi _{\mu }(a)\Omega _{\mu }
\end{equation*}
even in the case of unital $\ast $-algebras.

The dynamics $\alpha $ of a system $\mathbf{A}$ can be represented by a
unitary group $U$ on $H_{\mu }$ defined by extending
\begin{equation*}
U_{g}\gamma _{\mu }(a):=\gamma _{\mu }(\alpha _{g}(a)).
\end{equation*}
It satisfies
\begin{equation*}
U_{g}\pi _{\mu }(a)U_{g}^{\ast }=\pi _{\mu }(\alpha _{g}(a))
\end{equation*}
for all $g\in G$; also see \cite[Corollary 2.3.17]{BR}. The unitary
representation of $\beta $ will be denoted by $V$.

\begin{definition}
We call $\mathbf{F}$ a \emph{subsystem} of $\mathbf{A}$ if there exists an
injective unital $\ast$-homomorphism $\zeta$ of $F$ onto a von Neumann
subalgebra of $A$ such that $\mu\circ\zeta=\lambda$ and $\alpha_{g}\circ
\zeta=\zeta\circ\varphi_{g}$ for all $g\in G$. If $\zeta(F)$ is invariant
under $\sigma^{\mu}$, i.e. $\sigma_{t}^{\mu}(\zeta(F))=\zeta(F)$ for all
$t\in\mathbb{R}$, then $\mathbf{F}$ is called a \emph{modular subsystem} of
$\mathbf{A}$. If furthermore $\zeta:F\rightarrow A$ is surjective, then we
say that $\zeta$ is an \emph{isomorphism} of dynamical systems, and the
systems $\mathbf{A}$ and $\mathbf{F}$ are \emph{isomorphic}.
\end{definition}

The symbol $\zeta $ will be used uniformly for this purpose. We call it the
\emph{imbedding} of $\mathbf{F}$ into $\mathbf{A}$. Without loss one could
assume $\zeta $ to be an inclusion, so $F\subset A$, and we will often do
so. When the system $\mathbf{F}$ is a subsystem of $\mathbf{B}$, the
notation $\eta :F\rightarrow B$ instead of $\zeta $ will be used. If
$\mathbf{F}$ is a (modular) subsystem of both $\mathbf{A}$ and $\mathbf{B}$,
then we call it a \emph{common} (modular) subsystem of $\mathbf{A}$ and
$\mathbf{B}$.

Important nontrivial examples of modular subsystems will be studied in
Sections 4 and 5. In the meantime we note that in classical ergodic theory
all subsystems (also known as ``factors'') are of course modular. More
generally, if the state $\mu $ of our system $\mathbf{A}$ happens to be a
trace (i.e. $\mu (ab)=\mu (ba)$ for all $a,b\in A$), then again all
subsystems are modular. Also, in the physically relevant situation where the
dynamics $\alpha $ of $\mathbf{A}$ is the modular group of $\mu $, then any
subsystem of $\mathbf{A}$ is automatically modular.

A standard fact from Tomita-Takesaki theory related to modular subsystems is
the following (a proof of which is contained in \cite[Section 10.2]{St} for
example):

\begin{lemma}
Let $\mu$ be a faithful normal state on a von Neumann algebra $A$. Let $F$
be a von Neumann subalgebra of $A$, invariant under $\sigma^{\mu}$. Setting
$H_{F}:=\overline{\gamma_{\mu}(F)}$, we have $J_{\mu}H_{F}=H_{F}$.
\end{lemma}

When we construct relatively independent joinings of $\mathbf{A}$ and
$\mathbf{B}$ over a common modular subsystem $\mathbf{F}$ in the next
section, one of the key tricks (also used in \cite{Fid}) is to work with the
``commutant'' $\mathbf{\tilde{B}}$ of $\mathbf{B}$. Given a system $\mathbf{B}$,
this means the following: We set $\tilde{B}:=\pi _{\nu }(B)^{\prime }$
and then carry the state and dynamics of $\mathbf{B}$ over to $\tilde{B}$ in
a natural way using $j_{\nu }$ by defining a state $\tilde{\nu}$ and
$\ast$-automorphism $\tilde{\beta}_{g}$ on $\tilde{B}$ by $\tilde{\nu}(b):=\nu
\circ \pi _{\nu }^{-1}\circ j_{\nu }(b)$ and $\tilde{\beta}_{g}(b):=j_{\nu
}\circ \pi _{\nu }\circ \beta _{g}\circ \pi _{\nu }^{-1}\circ j_{\nu }(b)$
for all $g\in G$. From Tomita-Takesaki theory one has that $V_{g}J_{\nu
}=J_{\nu }V_{g}$ (see \cite[Construction 3.4]{D}). It follows that
\begin{equation*}
\tilde{\nu}(b)=\left\langle \Omega _{\nu },b\Omega _{\nu }\right\rangle
\end{equation*}
and
\begin{equation*}
\tilde{\beta}_{g}(b)=V_{g}bV_{g}^{\ast }
\end{equation*}
for all $b\in \tilde{B}$ and $g\in G$. In particular the latter tells us
that the unitary representation of $\tilde{\beta}$ is the same as that of
$\beta $, namely $V$.

\section{Relatively independent joinings}

Throughout this section we consider two systems $\mathbf{A}$ and $\mathbf{B}$
which have a common modular subsystem $\mathbf{F}$. We are going to
construct the relatively independent joining of $\mathbf{A}$ and $\mathbf{B}$
over $\mathbf{F}$. More precisely it will be a joining of $\mathbf{A}$ and
$\mathbf{\tilde{B}}$. Without loss we can assume $\zeta $ to be an inclusion,
so $F\subset A$.

Since $\mathbf{F}$ is a modular subsystem of $\mathbf{B}$, we obtain a
modular subsystem
$\mathbf{\tilde{F}}=\left(\tilde{F},\tilde{\lambda},\tilde{\varphi}\right)$
of $\mathbf{\tilde{B}}$ in the following way (as is easily checked): Set
\begin{equation*}
\tilde{F}:=j_{\nu}\circ\pi_{\nu}\circ\eta(F)\subset\tilde{B}
\end{equation*}
and let $\tilde{\lambda}:=\tilde{\nu}|_{\tilde{F}}$ and $\tilde{\varphi}_{g}
:=\tilde{\beta}_{g}|_{\tilde{F}}$ for all $g\in G$.

Since $F$ and $\tilde{F}$ are invariant under the modular groups of $\mu$
and $\tilde{\nu}$ respectively we know by Tomita-Takesaki theory (see for
example \cite[Theorem IX.4.2]{T}) that we have unique conditional
expectations
\begin{equation*}
D:A\rightarrow F
\end{equation*}
and
\begin{equation*}
\tilde{D}:\tilde{B}\rightarrow \tilde{F}
\end{equation*}%
such that $\lambda \circ D=\mu $ and $\tilde{\lambda}\circ \tilde{D}=\tilde{\nu}$.
(This corresponds to disintegrations in the classical case; see
\cite[Section 6.1]{G}.)

Setting $H_{\eta }:=\overline{\gamma _{\nu }(\eta (F))}$, we obtain a
well-defined unital $\ast $-homomorphism
\begin{equation*}
\pi _{\eta }:F\rightarrow B(H_{\eta }):a\mapsto \pi _{\nu }(\eta
(a))|_{H_{\eta }}
\end{equation*}
and one can check that $\left(H_{\eta},\pi_{\eta},\Omega_{\nu}\right)$ is a
cyclic representation of $\left( F,\lambda \right) $ and hence unitarily
equivalent \cite[Theorem 2.3.16]{BR} to $\left( H_{\lambda },\pi _{\lambda
},\Omega _{\lambda }\right) $. In other words there is a unique unitary
operator $u_{\eta }:H_{\lambda }\rightarrow H_{\eta }$ such that $u_{\eta
}\Omega _{\lambda }=\Omega _{\nu }$ and $u_{\eta }^{\ast }\pi
_{\eta}(a)u_{\eta }=\pi _{\lambda }(a)$ for all $a\in F$. From Lemma 2.4,
for any $a\in \tilde{F}$ we have $aH_{\eta }\subset H_{\eta }$, and
therefore we obtain a well-defined injective $\ast $-homomorphism
\begin{equation*}
\kappa _{\eta}:\tilde{F}\rightarrow B(H_{\lambda}):a\mapsto u_{\eta
}^{\ast}(a|_{H_{\eta}})u_{\eta}.
\end{equation*}
It is then easily shown from the definition of $\tilde{F}$ that $\kappa
_{\eta }(\tilde{F})\subset \pi _{\lambda }(F)^{\prime }$, and hence
\begin{equation*}
\delta :\pi _{\lambda }(F)\odot \kappa _{\eta }(\tilde{F})\rightarrow
B(H_{\lambda })
\end{equation*}
defined as the linear extension of
$\pi _{\lambda }(F)\times \kappa _{\eta }(\tilde{F})\rightarrow B(H_{\lambda }):(a,b)\mapsto ab$
is a well-defined unital $\ast $-homomorphism. We now introduce the
\emph{diagonal} state of $\lambda $ as the state
\begin{equation*}
\Delta _{\lambda }:F\odot \tilde{F}\rightarrow \mathbb{C}
\end{equation*}
defined by
\begin{equation*}
\Delta _{\lambda }(c):=\left\langle \Omega _{\lambda },\delta \circ \left(
\pi _{\lambda }\odot \kappa _{\eta }\right) (c)\Omega _{\lambda
}\right\rangle
\end{equation*}
for all $c\in F\odot \tilde{F}$. This enables us to define a linear
functional $\mu \odot _{\lambda }\tilde{\nu}$ on $A\odot \tilde{B}$ by
\begin{equation*}
\mu \odot _{\lambda }\tilde{\nu}:=\Delta _{\lambda }\circ (D\odot \tilde{D}).
\end{equation*}
(If we did not assume $\zeta $ to be an inclusion mapping, one would simply
replace $D$ by $\zeta ^{-1}\circ D$.)

This completes the basic construction. The next step is to show that $\mu
\odot _{\lambda }\tilde{\nu}$ is a joining of $\mathbf{A}$ and
$\mathbf{\tilde{B}}$. In fact we will show a bit more in the next proposition. We
first define a special class of joinings of $\mathbf{A}$ and $\mathbf{\tilde{B}}$:

\begin{definition}
Any $\omega\in J\left( \mathbf{A},\mathbf{\tilde{B}}\right) $ with
$\omega|_{F\odot\tilde{F}}=\Delta_{\lambda}$ is called a \emph{joining of}
$\mathbf{A}$ \emph{and} $\mathbf{\tilde{B}}$ \emph{over} $\mathbf{F}$. The
set of all such $\omega$ is denoted $J_{\lambda}\left( \mathbf{A},\mathbf{\tilde{B}}\right)$.
\end{definition}

\begin{proposition}
If $\mathbf{F}$ is a common modular subsystem of $\mathbf{A}$ and $\mathbf{B}
$, then $\mu\odot_{\lambda}\tilde{\nu}\in J_{\lambda}\left(\mathbf{A},\mathbf{\tilde{B}}\right)$.
\end{proposition}

\begin{proof}
From the uniqueness of $D$ it follows that $\varphi _{g}^{-1}\circ D\circ
\alpha _{g}=D$ and therefore $D\circ \alpha _{g}=\varphi _{g}\circ D$.
Similarly $\tilde{D}\circ \tilde{\beta}_{g}=\tilde{\varphi}_{g}\circ \tilde{D}$.
Denoting the unitary representation of $\varphi $ on $H_{\lambda }$ by
$W$, note that for $a\in F$ we have $u_{\eta }^{\ast }V_{g}|_{H_{\eta
}}u_{\eta }\pi _{\lambda }(a)u_{\eta }^{\ast }V_{g}|_{H_{\eta }}u_{\eta
}=u_{\eta }^{\ast }\pi _{\nu }(\beta _{g}(\eta (a)))u_{\eta }=\pi _{\lambda
}(\varphi _{g}(a))=W_{g}\pi _{\lambda }(a)W_{g}^{\ast }$. Letting this act
on the cyclic vector $\Omega _{\lambda }$, we obtain $u_{\eta }^{\ast
}V_{g}|_{H_{\eta }}u_{\eta }=W_{g}$. Therefore
$\kappa _{\eta }\left( \tilde{\varphi}_{g}(b)\right) =W_{g}\kappa _{\eta }(b)W_{g}^{\ast }$
for $b\in \tilde{F}$.

For $a\in A$ and $b\in \tilde{B}$ it follows that
\begin{eqnarray*}
&&\mu \odot _{\lambda }\tilde{\nu}\left( \alpha _{g}\odot \tilde{\beta}
_{g}(a\otimes b)\right) \\
&=&\left\langle \Omega _{\lambda },\pi _{\lambda }\left( \varphi
_{g}(D(a))\right) \kappa _{\eta }\left( \tilde{\varphi}_{g}(\tilde{D}
(b))\right) \Omega _{\lambda }\right\rangle \\
&=&\left\langle \Omega _{\lambda },W_{g}\pi _{\lambda }\left( D(a)\right)
W_{g}^{\ast }W_{g}\kappa _{\eta }\left( \tilde{D}(b)\right) W_{g}^{\ast
}\Omega _{\lambda }\right\rangle \\
&=&\mu \odot _{\lambda }\tilde{\nu}(a\otimes b)
\end{eqnarray*}
and therefore $\mu \odot _{\lambda }\tilde{\nu}\circ \left( \alpha _{g}\odot
\tilde{\beta}_{g}\right) =\mu \odot _{\lambda }\tilde{\nu}$ as required.
Since $D$ and $\tilde{D}$ are conditional expectations, it follows that
$D\odot \tilde{D}$ is positive (see \cite[p. 119]{St} for example). It is now
easily seen that $\mu \odot _{\lambda }\tilde{\nu}\in J\left( \mathbf{A},
\mathbf{\tilde{B}}\right) $. From the definition of $\mu \odot _{\lambda }
\tilde{\nu}$ it then immediately follows that
$\mu \odot _{\lambda }\tilde{\nu}\in J_{\lambda }\left( \mathbf{A},\mathbf{\tilde{B}}\right) $.
\end{proof}

\begin{definition}
The joining $\mu \odot _{\lambda }\tilde{\nu}$ of $\mathbf{A}$ and
$\mathbf{\tilde{B}}$ is called the \emph{relatively independent joining of}
$\mathbf{A}$ \emph{and} $\mathbf{\tilde{B}}$ \emph{over} $\mathbf{F}$. If $J_{\lambda
}\left( \mathbf{A},\mathbf{\tilde{B}}\right) =\left\{\mu\odot_{\lambda}
\tilde{\nu}\right\}$, then $\mathbf{A}$ and $\mathbf{\tilde{B}}$ are called
\emph{disjoint over} $\mathbf{F}$.
\end{definition}

In the remainder of this section we study simple but useful properties of
this relatively independent joining. These properties will be used in
subsequent sections of the paper.

\begin{proposition}
If $F=\mathbb{C}$, it follows that $\mu\odot_{\lambda}\tilde{\nu}=\mu \odot
\tilde{\nu}$ and $J_{\lambda}\left( \mathbf{A},\mathbf{\tilde{B}}\right)
=J\left( \mathbf{A},\mathbf{\tilde{B}}\right) $.
\end{proposition}

\begin{proof}
This follows directly from the definitions above.
\end{proof}

\begin{proposition}
Let $\mathbf{F}$ and $\mathbf{R}$ both be common modular subsystems of
$\mathbf{A}$ and $\mathbf{B}$, and $\mathbf{F}$ a subsystem of $\mathbf{R}$,
with inclusions $F\subset R\subset A$ giving the corresponding imbeddings of
the subsystems into $\mathbf{A}$, while $\theta :R\rightarrow B$ and $\eta
=\theta |_{F}$ give the imbeddings into $\mathbf{B}$. Then $J_{\psi }\left(
\mathbf{A},\mathbf{\tilde{B}}\right) \subset J_{\lambda }\left( \mathbf{A},
\mathbf{\tilde{B}}\right) $. Furthermore, if $\mu \odot _{\psi }\tilde{\nu}
|_{R\odot \tilde{R}}=\mu \odot _{\lambda }\tilde{\nu}|_{R\odot \tilde{R}}$,
then $R=F$.
\end{proposition}

\begin{proof}
For the first part we only need to show that $\Delta _{\psi }|_{F\odot
\tilde{F}}=\Delta _{\lambda }$. Given $\left( H_{\psi },\pi _{\psi },\Omega
_{\psi }\right) $, we can assume that $\left( H_{\lambda },\pi _{\lambda
},\Omega _{\lambda }\right) $ is given by $H_{\lambda }=\overline{\gamma
_{\psi }(F)}$, $\pi _{\lambda }(a)=\pi _{\psi }(a)|_{H_{\lambda }}$ for
$a\in F$, and $\Omega _{\lambda }=\Omega _{\psi }$, without changing $\Delta
_{\lambda }$, because of unitary equivalence between different cyclic
representations of the same state. We can also define $\kappa _{\theta }:
\tilde{R}\rightarrow B(H_{\psi })$ in the same way as $\kappa _{\eta }$
above, with the corresponding unitary operator denoted by $u_{\theta }$
instead of $u_{\eta }$. Note that $u_{\theta }\gamma _{\lambda }(a)=\pi
_{\theta }(a)\Omega _{\nu }=\pi _{\eta }(a)\Omega _{\nu }=u_{\eta }\gamma
_{\lambda }(a)$ for $a\in F$, from which follows that $\kappa _{\theta
}(b)|_{H_{\lambda }}=\kappa _{\eta }(b)$ for $b\in \tilde{F}$,\ hence
$\Delta _{\psi }|_{F\odot \tilde{F}}=\Delta _{\lambda }$.

Assuming $\mu \odot _{\psi }\tilde{\nu}|_{R\odot \tilde{R}}=\mu \odot
_{\lambda }\tilde{\nu}|_{R\odot \tilde{R}}$, we have (in $H_{\psi }$)
\begin{align*}
\left\langle \pi _{\psi }(a)\Omega _{\psi },\kappa _{\theta }(b)\Omega
_{\psi }\right\rangle & =\mu \odot _{\psi }\tilde{\nu}(a^{\ast }\otimes
b)=\mu \odot _{\lambda }\tilde{\nu}(a^{\ast }\otimes b) \\
& =\mu \odot _{\lambda }\tilde{\nu}(a^{\ast }\otimes \tilde{D}(b))=\mu \odot
_{\psi }\tilde{\nu}(a^{\ast }\otimes \tilde{D}(b)) \\
& =\left\langle \pi _{\psi }(a)\Omega _{\psi },\kappa _{\theta }(\tilde{D}
(b))\Omega _{\psi }\right\rangle
\end{align*}
for any $a\in R$ and $b\in \tilde{R}$ (with $\tilde{R}$ of course defined in
an analogous way to $\tilde{F}$). Since $\Omega _{\psi }$ is cyclic for $\pi
_{\psi }$, and therefore separating for $\kappa _{\theta }(\tilde{R})\subset
\pi _{\psi }(R)^{\prime }$, it follows that $\kappa _{\theta }(b)=\kappa
_{\theta }(\tilde{D}(b))$. Therefore $b|_{H_{\theta }}=\tilde{D}%
(b)|_{H_{\theta }}$ where $H_{\theta }:=\overline{\gamma _{\nu }(\theta (R))}
$, and hence $b\Omega _{\nu }=\tilde{D}(b)\Omega _{\nu }$. Since $\Omega
_{\nu }$ is separating for $\tilde{B}$ we conclude that $b=\tilde{D}(b)\in
\tilde{F}$ and therefore $\tilde{R}=\tilde{F}$. It follows that $R=F$ as
required.
\end{proof}

Finally we consider a Hilbert space characterization of relatively
independent joinings. We will work in the following setting: Let $\omega $
be a coupling of $\left( A,\mu \right) $ and $\left( \tilde{B},\tilde{\nu}
\right) $ as defined in Section 2. Then we can ``imbed'' $\left( H_{\mu
},\pi _{\mu },\Omega _{\mu }\right) $ and
$\left( H_{\tilde{\nu}},\pi _{\tilde{\nu}},\Omega _{\tilde{\nu}}\right) $
into $\left( H_{\omega },\pi_{\omega},\Omega _{\omega }\right) $ in a natural
way (see \cite[Construction 2.3]{D} for explicit details; also note that
$\left( H_{\tilde{\nu}},\pi _{\tilde{\nu}},\Omega _{\tilde{\nu}}\right) $ is unitarily
equivalent to $\left( H_{\nu},\iota _{\tilde{B}},\Omega _{\nu }\right) $ but
not necessarily equal). In particular $H_{\mu }$ and $H_{\tilde{\nu}}$ are
then subspaces of $H_{\omega }$ and the corresponding cyclic vectors are
equal to $\Omega _{\omega }$. We set $H_{\lambda }^{\prime }=\overline{
\gamma _{\mu }(F)}\subset H_{\mu }$ and $H_{\tilde{\lambda}}=\overline{
\gamma _{\tilde{\nu}}(\tilde{F})}\subset H_{\tilde{\nu}}$,

We are particularly interested in the case where $\omega |_{F\odot \tilde{F}
}=\Delta _{\lambda }$. If this condition is satisfied one has $H_{\lambda
}^{\prime }=H_{\tilde{\lambda}}$ as follows: $H_{\lambda }^{\prime }$ and
$H_{\tilde{\lambda}}$ are both contained $H_{\odot }:=\overline{\gamma
_{\omega }\left( F\odot \tilde{F}\right) }$. Consider $x\in H_{\odot
}\ominus H_{\lambda }^{\prime }$ and a sequence $b_{n}\in F\odot \tilde{F}$
such that $\gamma _{\omega }(b_{n})\rightarrow x$, then for every $a\in F$
it follows from $\omega |_{F\odot \tilde{F}}=\Delta _{\lambda }$ that $
0=\left\langle \gamma _{\mu }(a),x\right\rangle =\left\langle \gamma
_{\lambda }(a),x^{\prime }\right\rangle $ where the limit $x^{\prime
}:=\lim_{n\rightarrow \infty }\delta \circ \left( \pi _{\lambda }\odot
\kappa _{\eta }\right) (b_{n})\Omega _{\lambda }$ exists in $H_{\lambda }$
by unitary equivalence between the cyclic representations $\left( H_{\lambda
},\delta \circ \left( \pi _{\lambda }\odot \kappa _{\eta }\right) ,\Omega
_{\lambda }\right) $ and $\left( H_{\odot },\pi _{\omega }(\cdot
)|_{H_{\odot }},\Omega _{\omega }\right) $ of $\left( F\odot \tilde{F}
,\Delta _{\lambda }\right) $. Therefore $x^{\prime }=0$, and so $x=0$ by the
same unitary equivalence. Hence $H_{\lambda }^{\prime }=H_{\odot }$. In a
similar fashion $H_{\tilde{\lambda}}=H_{\odot }$ proving the claim. By
unitary equivalence one can choose $\left( H_{\lambda },\pi _{\lambda
},\Omega _{\lambda }\right) $ to be $\left( H_{\lambda }^{\prime },\pi _{\mu
}(\cdot )|_{H_{\lambda }^{\prime }},\Omega _{\mu }\right) $, with a
corresponding change in $\kappa _{\eta }$, without changing the state
$\Delta _{\lambda }$ . In conclusion we therefore have
\begin{equation*}
H_{\tilde{\lambda}}=H_{\lambda }
\end{equation*}
when $\omega |_{F\odot \tilde{F}}=\Delta _{\lambda }$.

In this setting we now have the following result:

\begin{proposition}
Suppose that $\omega$ is a coupling of $\left( A,\mu\right) $ and
$\left(\tilde{B},\tilde{\nu}\right) $ such that $\omega|_{F\odot\tilde{F}
}=\Delta_{\lambda}$. Then $\omega=\mu\odot_{\lambda}\tilde{\nu}$ if and only
if any of the following three equivalent conditions are satisfied: $\left(
H_{\mu}\ominus H_{\lambda}\right) \perp\left( H_{\tilde{\nu}}\ominus
H_{\lambda}\right) $, $\left( H_{\mu}\ominus H_{\lambda}\right) \perp H_{\tilde{\nu}}$,
or $H_{\mu}\perp\left(H_{\tilde{\nu}}\ominus H_{\lambda}\right)$.
\end{proposition}

\begin{proof}
The equivalence of the three conditions is easily verified. So assume
$\omega =\mu \odot _{\lambda }\tilde{\nu}$ and consider any $x\in H_{\mu
}\ominus H_{\lambda }$ and $y\in H_{\tilde{\nu}}$, as well as sequences
$a_{n}\in A$ and $b_{n}\in \tilde{B}$ such that $\gamma _{\mu }\left(
a_{n}\right) \rightarrow x$ and $\gamma _{\tilde{\nu}}\left( b_{n}\right)
\rightarrow y$. Then
\begin{align*}
\left\langle x,\gamma _{\tilde{\nu}}(b_{n})\right\rangle &
=\lim_{m\rightarrow \infty }\omega \left( a_{m}^{\ast }\otimes b_{n}\right)
=\lim_{m\rightarrow \infty }\mu \odot _{\lambda }\tilde{\nu}\left(
a_{m}^{\ast }\otimes b_{n}\right) \\
& =\left\langle \Omega _{\lambda },\left( \pi _{\lambda }\circ D(a_{m}^{\ast
})\right) \left( \kappa _{\eta }\circ \tilde{D}(b_{n})\right) \Omega
_{\lambda }\right\rangle \\
& =\left\langle \Omega _{\lambda },\left( \pi _{\lambda }\circ D(a_{m}^{\ast
})\right) \left( \kappa _{\eta }\circ \tilde{D}\left( \tilde{D}
(b_{n})\right) \right) \Omega _{\lambda }\right\rangle \\
& =\left\langle x,\gamma _{\tilde{\nu}}\left( \tilde{D}(b_{n})\right)
\right\rangle \\
& =0
\end{align*}
and therefore $\left\langle x,y\right\rangle =0$. It follows that $\left(
H_{\mu }\ominus H_{\lambda }\right) \perp H_{\tilde{\nu}}$. Conversely,
assume that
$\left( H_{\mu }\ominus H_{\lambda }\right) \perp \left( H_{\tilde{\nu}}\ominus H_{\lambda }\right) $
and consider any $a\in A$ and $b\in \tilde{B}$. Let $P$ and $\tilde{P}$
respectively be the projections of $H_{\mu }$ and
$H_{\tilde{\nu}}$ on $H_{\lambda }$. Then for any $a\in A$
and $b\in \tilde{B}$
\begin{align*}
\omega \left( a\otimes b\right) & =\left\langle \gamma _{\mu }(a^{\ast}),
\gamma _{\tilde{\nu}}(b)\right\rangle =\left\langle P\gamma _{\mu}(a^{\ast }),
\tilde{P}\gamma _{\tilde{\nu}}(b)\right\rangle \\
& =\left\langle \gamma _{\mu }\left( D(a^{\ast })\right) ,
\gamma _{\tilde{\nu }}\left( \tilde{D}(b)\right) \right\rangle =
\omega \left( D(a)\otimes
\tilde{D}(b)\right) \\
& =\Delta _{\lambda }\left( D(a)\otimes \tilde{D}(b)\right) =\mu \odot
_{\lambda }\tilde{\nu}\left( a\otimes b\right)
\end{align*}
which is sufficient.
\end{proof}

We will also unitarily represent the dynamics of $\mathbf{\tilde{B}}$ on
$H_{\tilde{\nu}}$, and we denote this representation by $\tilde{V}$. This
representation is used in Sections 4 and 5 instead of $V$ on $H_{\nu }$ to
fit into the setting of Proposition 3.6.

\section{Relative ergodicity}

In this section and the next we study how the relatively independent
joinings constructed in the previous section relate to properties of
systems. In particular in this section we consider relative ergodicity,
which is a simple generalization of ergodicity. See for example
\cite[Section 6.6]{G} for a discussion of the classical case.

The relevant terminology and definitions are as follows: The \emph{fixed
point algebra} of a system $\mathbf{A}$ (or, put differently, of $\alpha $)
is defined as
\begin{equation*}
A^{\alpha }:=\left\{ a\in A:\alpha _{g}(a)=a\text{ for all }g\in G\right\}
\end{equation*}
and this gives an identity subsystem $\mathbf{A}^{\alpha }$ of $\mathbf{A}$
by simply restricting the state of $\mathbf{A}$ to $A^{\alpha }$ . We call
$\mathbf{A}^{\alpha }$ the \emph{fixed point subsystem} of $\mathbf{A}$.
Similarly, we call
\begin{equation*}
H_{\mu }^{U}:=\left\{ x\in H_{\mu }:U_{g}x=x\text{ for all }g\in G\right\}
\end{equation*}
the fixed point space of $U$.

\begin{definition}
Let $\mathbf{F}$ be a subsystem of $\mathbf{A}$. We call $\mathbf{A}$
\emph{ergodic relative to} $\mathbf{F}$ if $A^{\alpha}\subset\zeta(F)$.
\end{definition}

If $\mathbf{F}$ in Definition 4.1 is the trivial system, then $\mathbf{A}$
is simply called \emph{ergodic}.

We note the following important facts:

\begin{proposition}
The fixed point subsystem $\mathbf{A}^{\alpha}$ is a modular subsystem of
$\mathbf{A}$. Furthermore, $H_{\mu}^{U}=\overline{\gamma_{\mu}(A^{\alpha})}$.
\end{proposition}

\begin{proof}
By \cite[Corollary VIII.1.4]{T} $\alpha_{g}\circ\sigma_{t}^{\mu}=\alpha
_{g}\circ\sigma_{t}^{\mu\circ\alpha_{g}}=\sigma_{t}^{\mu}\circ\alpha_{g}$
from which $\sigma_{t}^{\mu}(A^{\alpha})=A^{\alpha}$ follows.

We now show that $H_{\mu }^{U}=\overline{\gamma _{\mu }(A^{\alpha })}$. It
is clear that $\overline{\gamma _{\mu }(A^{\alpha })}\subset H_{\mu }^{U}$.
The converse follows from the Kov\'{a}cs-Sz\"{u}cs mean ergodic theorem
\cite[Proposition 4.3.8]{BR} (also see \cite{KS}) using a argument similar
to that of \cite[Theorem 4.3.20]{BR} (also see \cite{Jad}): Let $P$ be the
projection of $H_{\mu }$ onto $H_{\mu }^{U}$, then, since $\Omega _{\mu }\in
H_{\mu }^{U}$ is cyclic for $\pi _{\mu }(A)^{\prime }$, there exists a
unique normal $(\pi _{\mu }\circ \alpha _{g}\circ \pi _{\mu }^{-1})_{g\in G}$
invariant projection $E:\pi _{\mu }(A)\rightarrow \pi _{\mu }(A^{\alpha })$.
Furthermore, $E$ has the property that $Ea$ is the unique element of $\pi
_{\mu }(A)$ such that $\left( Ea\right) P=PaP$, for every $a\in \pi _{\mu
}(A)$. Hence $(Ea)\Omega _{\mu }=Pa\Omega _{\mu }$ from which we obtain
$\overline{\gamma _{\mu }(A^{\alpha })}\supset H_{\mu }^{U}$.
\end{proof}

Now for the main result of this section:

\begin{theorem}
Let the identity system $\mathbf{F}$ be a modular subsystem of $\mathbf{A}$.
If $\mathbf{A}$ and $\mathbf{\tilde{B}}$ are disjoint over $\mathbf{F}$ for
$\mathbf{B}=\mathbf{A}^{\alpha }$, then $\mathbf{A}$ is ergodic relative to
$\mathbf{F}$. On the other hand, if $\mathbf{A}$ is ergodic relative to
$\mathbf{F}$, then $\mathbf{A}$ and $\mathbf{\tilde{B}}$ are disjoint over
$\mathbf{F}$ for every identity system $\mathbf{B}$ which has $\mathbf{F}$ as
a modular subsystem.
\end{theorem}

\begin{proof}
Without loss we assume that $\zeta $ is an inclusion, so $F\subset A$. Now,
suppose $\mathbf{A}$ is not ergodic relative to $\mathbf{F}$, in other words
$F$ is strictly contained in $A^{\alpha }$. Note that $\mathbf{F}$ is a
modular subsystem of $\mathbf{A}^{\alpha }$, since $\mathbf{F}$ and
$\mathbf{A}^{\alpha }$ are modular subsystems of $\mathbf{A}$. Then apply
Proposition 3.5 with $\mathbf{B}=\mathbf{R}=\mathbf{A}^{\alpha }$ to obtain
$\mu \odot _{\psi }\tilde{\nu}\neq \mu \odot _{\lambda }\tilde{\nu}$. But
both these joinings are contained in
$J_{\lambda }\left( \mathbf{A},\mathbf{\tilde{B}} \right) $
by Propositions 3.2 and 3.5, so $\mathbf{A}$ and $\mathbf{\tilde{B}}$
are not disjoint over $\mathbf{F}$.

Conversely, assume that $\mathbf{A}$ is ergodic relative to $\mathbf{F}$,
and let $\mathbf{B}$ be any identity system which has $\mathbf{F}$ as a
modular subsystem. We are going to apply Proposition 3.6. So consider any
$\omega \in J_{\lambda }(\mathbf{A},\mathbf{\tilde{B}})$. From this joining
we obtain a conditional expectation operator
$P_{\omega }:H_{\tilde{\nu}}\rightarrow H_{\mu }$
(i.e. $\left\langle x,P_{\omega }y\right\rangle
=\left\langle x,y\right\rangle $ for all $x\in H_{\mu }$ and
$y\in H_{\tilde{\nu}}$) such that $U_{g}P_{\omega }=P_{\omega }\tilde{V}_{g}$
(see \cite[Proposition 2.4]{D}). Then for any $b\in \tilde{B}$ we have
$U_{g}P_{\omega }\gamma _{\tilde{\nu}}(b)=P_{\omega }\gamma _{\tilde{\nu}}(b)$
, since $\mathbf{\tilde{B}}$ is an identity system. Since $F=A^{\alpha }$,
it follows from Proposition 4.2 that $P_{\omega }\gamma _{\tilde{\nu}}(b)\in
\overline{\gamma _{\mu }(A^{\alpha })}=H_{\lambda }$. Noting that $P_{\omega
}$ is the projection of $H_{\omega }$ onto $H_{\mu }$ restricted to
$H_{\tilde{\nu}}$, we conclude that $H_{\tilde{\nu}}\perp \left( H_{\mu }\ominus
H_{\lambda }\right) $, and hence $\omega =\mu \odot _{\lambda }\tilde{\nu}$
by Proposition 3.6.
\end{proof}

This result contains \cite[Theorem 2.1]{D2} as a special case by Proposition
3.4, namely $\mathbf{A}$ is ergodic if and only if it is disjoint (over the
trivial system) from all identity systems. Also note that in Theorem 4.3,
$\mathbf{A}$ being ergodic relative to $\mathbf{F}$ means exactly that
$A^{\alpha }=\zeta (F)$, i.e. $\mathbf{F}$ is isomorphic to $\mathbf{A}
^{\alpha }$, since $\mathbf{F}$ is assumed to be an identity system.

\section{Compact subsystems}

In this section we study an analogue of the previous section in terms of
compact systems instead of identity systems. We therefore proceed by first
building up some background regarding compact subsystems, which includes
extending certain results from \cite[Section 4]{NSZ} to more general group
actions in a way (see Lemma 5.3 and Theorem 5.4 below) that is relevant for
our intended applications. Our proofs of these extensions follow the basic
plans of those in \cite{NSZ}, but many of the details differ. (To avoid any
confusion, we note that our terminology regarding topologies on $B(H)$
differs slightly from that of \cite{NSZ} and furthermore we do not always
use the same topologies that \cite{NSZ} used in their results; we use the
terminology of \cite[Section 2.4.1]{BR}, in particular the \emph{weak}
topology on $B(H)$ is generated by $\left\langle x,(\cdot )y\right\rangle $
for $x,y\in H$.) The additional definitions that we use in this section are
the following:

\begin{definition}
A system $\mathbf{A}$ is \emph{compact} if the orbit $U_{G}x$ is totally
bounded (i.e. $\overline{U_{G}x}$ is compact) in $H_{\mu }$ for every $x\in
H_{\mu }$. An \emph{eigenvector} of $U$ is an $x\in H_{\mu }\backslash \{0\}$
such that there is a function, called its \emph{eigenvalue}, $\chi
:G\rightarrow \mathbb{C}$ such that $U_{g}x=\chi (g)x$ for all $g\in G$. The
set of all eigenvalues is denoted by $\sigma _{\mathbf{A}}$. Let $H_{0}$
denote the Hilbert subspace of $H_{\mu }$ spanned by the eigenvectors of $U$.
We will call $u\in A\backslash \{0\}$ an \emph{eigenoperator} of $\alpha $
if there is a function $\chi :G\rightarrow \mathbb{C}$, its \emph{eigenvalue},
such that $\alpha _{g}(u)=\chi (g)u$ for all $g\in G$.
\end{definition}

As already mentioned, our immediate goal is to study compact subsystems of a
given system. We start with an analogue of the first part of Proposition 4.2:

\begin{proposition}
For a system $\mathbf{A}$, denote by $A^{K}$ the von Neumann subalgebra of
$A$ generated by the eigenoperators of $\alpha $. Then $\alpha
_{g}(A^{K})=A^{K}$, which allows us to define a subsystem
$\mathbf{A}^{K}=\left( A^{K},\mu ^{K},\alpha ^{K}\right) $ of $\mathbf{A}$
(in terms of an inclusion as the imbedding into $\mathbf{A}$) by setting
$\mu ^{K}:=\mu |_{A^{K}}$ and $\alpha _{g}^{K}(a):=\alpha _{g}(a)$ for all
$a\in A^{K}$. The system $\mathbf{A}^{K}$ is a compact modular subsystem of
$\mathbf{A}$.
\end{proposition}

\begin{proof}
Let $S$ be the $\ast $-algebra generated by the eigenoperators of $\alpha $.
Then clearly $\alpha _{g}(S)=S\subset A^{K}$. But $A^{K}=S^{\prime \prime }$
, so $S$ is $\sigma $-weakly dense in $A^{K}$, while $\alpha _{g}$ is
$\sigma $-weakly continuous for each $g$. It follows that
$\alpha_{g}(A^{K})=A^{K}$.

Next we prove that the subsystem $\mathbf{A}^{K}$ so obtained, is compact.
Note that for any eigenoperator $u$ of $\alpha $, $\gamma _{\mu }(u)$ is an
eigenvector of $U$. Furthermore, it is easily seen that if $u$ and $v$ are
eigenoperators of $\alpha $, then $u^{\ast }$ is also eigenoperator of
$\alpha $, while $uv$ is either zero or an eigenoperator of $\alpha $. It
follows that $\gamma _{\mu }(S)\subset H_{0}$, hence $\gamma _{\mu
}(A^{K})\subset H_{0}$. It is easy to show from the definition of $H_{0}$
that each of its elements has a totally bounded (i.e. relatively compact)
orbit under $U$. Also note that $\left( A^{K},\mu \right) $ can be
cyclically represented on $H_{\mu }^{K}:=\overline{\gamma _{\mu }(A^{K})}$,
which is contained in $H_{0}$, with the unitary representation of $\alpha
^{K}$ given by the restriction of $U_{g}$ to $H_{\mu }^{K}$. Therefore
$\mathbf{A}^{K}$ is indeed compact.

Lastly note that as in Proposition 4.2,\ for any eigenoperator $u$ of
$\alpha $ with eigenvalue $\chi $, we have that $\theta _{g}(\sigma _{t}^{\mu
}(u))=\sigma _{t}^{\mu }(\theta _{g}(u))=\chi (g)\sigma _{t}^{\mu }(u)$ so
$\sigma _{t}^{\mu }(u)\in S$. It follows that $\sigma _{t}^{\mu }(S)=S$,
hence $\sigma _{t}^{\mu }(A^{K})=A^{K}$, i.e. $\mathbf{A}^{K}$ is a modular
subsystem of $\mathbf{A}$.
\end{proof}

In order to go further, we first prove a technical result which is a version
of \cite[Lemma 4.1]{NSZ} appropriate for our needs:

\begin{lemma}
Let $H$ be a Hilbert space and $T$ a weakly closed vector subspace of $B(H)$.
Assume that there is a unit vector $\Omega \in H$ such that
$\overline{T^{\prime}\Omega }=H$. Let $G$ be an arbitrary group and $W_{g}:H\rightarrow
H$ a linear isometry for every $g\in G$ such that $W_{g}W_{h}=W_{gh}$ for
all $g,h\in G$, and with
\begin{equation*}
K_{a}:=\overline{\left\{ W_{g}a\Omega :g\in G\right\} }
\end{equation*}
compact in $H$ for every $a\in T$. Endow $B(T)$ with the topology of
pointwise strong convergence generated by the seminorms $\Theta \mapsto
\left\| \Theta (a)x\right\| $ where $a\in T$ and $x\in H$. Let $\Gamma$ be
the set of all linear contractions $\Theta \in B(T)$ such that $\Theta
(a)\Omega \in K_{a}$ for all $a\in T$, and view it as a topological subspace
of $B(T)$. Then $\Gamma $ is a compact group with composition as its
operation and the identity map on $T$ as its identity.
\end{lemma}

\begin{proof}
Let us first consider the topology on $\Gamma $. For $a\in T,a^{\prime }\in
T^{\prime },x\in H$ and $\Theta _{1},\Theta _{2}\in B(T)$ we have
\begin{equation*}
\left\| \Theta _{1}(a)x-\Theta _{2}(a)x\right\| \leq 2\left\| a\right\|
\left\| x-a^{\prime }\Omega \right\| +\left\| a^{\prime }\right\| \left\|
\Theta _{1}(a)\Omega -\Theta _{2}(a)\Omega \right\|
\end{equation*}
from which (together with $\overline{T^{\prime }\Omega }=H$) it follows that
the topology on $B(T)$ is in fact generated by the seminorms $\Theta \mapsto
\left\| \Theta (a)\Omega \right\| $ where $a\in T$, hence we have a
simplified description of the topology on $\Gamma $. Now consider the
compact space
\begin{equation*}
K:=\prod_{a\in T}K_{a}
\end{equation*}
and the function
\begin{equation*}
f:\Gamma \rightarrow K:\Theta \mapsto \left( \Theta (a)\Omega \right) _{a\in
T}
\end{equation*}
which is easily shown to be injective because of $\overline{T^{\prime
}\Omega }=H$. From the simplified description of the topology on $\Gamma $
described above and the definition of the product topology on $K$ it is
clear that $f:\Gamma \rightarrow f(\Gamma )$ is a homeomorphism. We now show
that $f(\Gamma )$ is closed:

Let $x=\left( x(a)\right) _{a\in T}$ be in the closure of $f(\Gamma
)$. Let $\left( \Theta _{\iota }\right) _{\iota }$ be a net in
$\Gamma $ such that $\left( f(\Theta _{\iota })\right) _{\iota }$
converges to $x$, then from the inequality above $\left( \Theta
_{\iota }(a)y\right) _{\iota }$ is a Cauchy net in $H$ for every
$y\in H$ and $a\in T$. This allows us to define $\Theta
(a)y:=\lim_{\iota }\Theta _{\iota }(a)y$ and since each $\Theta
_{\iota }$ is a contraction, $\Theta (a)\in B(H)$ and $\left\|
\Theta (a)\right\| \leq \left\| a\right\| $ . By this definition
$\left( \Theta _{\iota }(a)\right) _{\iota }$ converges strongly and
therefore also weakly to $\Theta (a)$, but $T$ is weakly closed
therefore $\Theta (a)\in T$. Since $\Theta _{\iota }(a)\Omega \in
K_{a}$ and $K_{a}$ is closed, it again follows from the definition
of $\Theta (a)$ that $\Theta (a)\Omega \in K_{a}$ for all $a\in T$.
This shows that $\Theta \in \Gamma $. But by our choice of $\left(
\Theta _{\iota }\right) _{\iota }$ we have $x(a)=\lim_{\iota }\Theta
_{\iota }(a)\Omega =\Theta (a)\Omega $ so $x=f(\Theta )\in
f(\Gamma)$ and therefore $f(\Gamma )$ is indeed closed.

It follows that $f(\Gamma)$ is compact, since $K$ is, and therefore $\Gamma$
is compact. Also note that $\Gamma$ is Hausdorff.

Now we start to prove the group structure on $\Gamma $. Note that $W_{1}=1$,
since $W_{g}W_{1}=W_{g}$ while $W_{g}$ is injective. It follows that
id$_{T}\in \Gamma $. For any $\Theta _{1},\Theta _{2}\in \Gamma $ we have
\begin{equation*}
\Theta _{1}\circ \Theta _{2}(a)\Omega \in K_{\Theta _{2}(a)}\subset
\overline{\left\{ W_{g}W_{h}a\Omega :g,h\in G\right\} }=K_{a}
\end{equation*}
and therefore $\Theta _{1}\circ \Theta _{2}\in \Gamma $.

Next we show the continuity of this product on $\Gamma$. For
$\Theta\in\Gamma $ we have $\Theta(a)\Omega\in K_{a}\subset\left\{x\in
H:\left\| x\right\| =\left\| a\Omega\right\|\right\}$, since each $W_{g}$ is
an isometry. Therefore we have a well-defined linear isometry
\begin{equation*}
W_{\Theta}:\overline{T\Omega}\rightarrow\overline{T\Omega}
\end{equation*}
such that $W_{\Theta}(a\Omega)=\Theta(a)\Omega$, so in effect $W_{\Theta}$
represents $\Theta$ on $\overline{T\Theta}$. Using $W_{\Theta}$ one can show
that
\begin{align*}
& \left\| (\Theta_{1}\circ\Theta_{1}^{\prime})(a)\Omega-
(\Theta_{2}\circ\Theta_{2}^{\prime})(a)\Omega\right\| \\
& \leq\left\|
\Theta_{1}^{\prime}(a)\Omega-\Theta_{2}^{\prime}(a)\Omega\right\| +\left\|
\Theta_{1}\left( \Theta_{2}^{\prime}(a)\right) \Omega-
\Theta_{2}\left(\Theta_{2}^{\prime}(a)\right) \Omega\right\|
\end{align*}
for all $\Theta_{1},\Theta_{2},\Theta_{1}^{\prime},\Theta_{2}^{\prime}\in
\Gamma$ and $a\in T$. From this it follows that the product on $\Gamma$ is
indeed continuous.

We now complete the proof that $\Gamma $ is a group by considering inverses.
We have already shown that $\Gamma $ is a semigroup, so for $\Theta \in
\Gamma $ it follows that $\Theta ^{n}\in \Gamma $ for every $n\in \mathbb{N}
=\left\{ 1,2,3,...\right\}$. Furthermore $W_{\Theta }^{n}(a\Omega )=\Theta
^{n}(a)\Omega \in K_{a}$, so since $K_{a}$ is compact, the orbit
$\left(W_{\Theta}^{n}(a\Omega)\right)_{n\in \mathbb{N}}$ is relatively
compact (and therefore totally bounded) for every $a\in T$. Since $W_{\Theta
}$ is an isometry, it follows from \cite[Corollary 9.10]{NSZ} (and remarks
made just after it) that for every $\varepsilon >0$ and every finite set
$E\subset T$ there exists an $n(E,\varepsilon )\in \mathbb{N}$ such that
\begin{equation*}
\left\| \Theta ^{n(E,\varepsilon )}(a)\Omega -a\Omega \right\| <\varepsilon
\end{equation*}
for all $a\in E$. Defining $(E,\varepsilon)\leq(E^{\prime},\varepsilon
^{\prime })$ to mean $E\subset E^{\prime }$ and $\varepsilon \geq\varepsilon
^{\prime }$, we obtain a net $\left( \Theta ^{n(E,\varepsilon )}\right)
_{(E,\varepsilon )}$ in $\Gamma $ which is seen to converge to id$_{T}$ by
using the description of the topology given at the beginning of this proof.
Since $\Gamma $ is compact, the net $\left(\Theta^{n(E,\varepsilon
)-1}\right)_{(E,\varepsilon )}$ has a limit point in $\Gamma $, say $\Theta
^{\prime }$, and it is then not too difficult to show that $\Theta \circ
\Theta ^{\prime }=$ id$_{T}=\Theta ^{\prime }\circ \Theta $. In other words
every element of $\Gamma $ has an inverse in $\Gamma $, and therefore
$\Gamma $ is a group.

By the single theorem appearing in \cite{E}, it follows that $\Gamma$ is a
topological group, which completes the proof.
\end{proof}

Now we can state and prove the basic result regarding compact subsystems,
which is the version of \cite[Theorem 4.2]{NSZ} that we will use. In this
result and the rest of the section our group $G$ is assumed to be abelian,
although this assumption is unnecessary in the first part of the proof of
the next theorem, as will be indicated.

\begin{theorem}
Consider a system $\mathbf{A}$ with $G$ assumed to be abelian. Set
\begin{align*}
T& :=\left\{ a\in \pi _{\mu }(A):a\Omega _{\mu }\in H_{0}\right\} \\
\xi (a)& :=\left\langle \Omega _{\mu },a\Omega _{\mu }\right\rangle \\
\tau _{g}(a)& :=U_{g}aU_{g}^{\ast }
\end{align*}%
for all $a\in T$ and $g\in G$. Then $\mathbf{T:=}\left( T,\xi ,\tau \right) $
is a subsystem of $\mathbf{A}$ which is isomorphic to $\mathbf{A}^{K}$ with
the isomorphism given by $\pi _{\mu }^{-1}|_{T}:T\rightarrow A^{K}$. It
follows that $\mathbf{A}^{K}$ is the largest compact subsystem of $\mathbf{A}
$ in the sense that $A^{K}$ contains the image of the algebra of any compact
subsystem $\mathbf{F}$ of $\mathbf{A}$ (under the imbedding of $\mathbf{F}$
in $\mathbf{A}$). Furthermore, $\overline{T\Omega _{\mu }}=H_{0}$, so
$\mathbf{T}$ can be cyclically represented on $H_{0}$ and
$\sigma _{\mathbf{T}}=\sigma _{\mathbf{A}}$.
\end{theorem}

\begin{proof}
The notation $M=\pi _{\mu }(A)$, $\theta _{g}:M\rightarrow M:a\mapsto
U_{g}aU_{g}^{\ast }$, $H=H_{\mu }$ and $\Omega =\Omega _{\mu }$ will be used
in what follows. We divide the proof into a number of parts.

(i) First we show (see also \cite[Proposition 3.2]{NSZ}) that for any $\chi
\in \mathbb{\sigma }_{\mathbf{A}}$ we have
\begin{equation*}
\overline{M_{\chi }\Omega }=H_{\chi }
\end{equation*}
where \
\begin{equation*}
M_{\chi }:=\left\{ a\in M:\theta _{g}(a)=\chi (g)a\text{ for all }g\in
G\right\}
\end{equation*}
and
\begin{equation*}
H_{\chi }:=\left\{ x\in H:U_{g}x=\chi (g)x\text{ for all }g\in G\right\} .
\end{equation*}
To do this, let $P_{\chi }$ be the projection of $H$ on $H_{\chi }$. Note
that $\left| \chi (g)\right| =1$ and $\chi (g)\chi (h)=\chi (gh)$. It
follows that $g\mapsto \overline{\chi (g)}U_{g}$ is also a unitary group on
$H$, and therefore by the Alaoglu-Birkhoff mean ergodic theorem
\cite[Proposition 4.3.4]{BR} (also see \cite{AB}) we know that $P_{\chi }$
is in the strong closure of the convex hull of $\left\{ \overline{\chi (g)}
U_{g}:g\in G\right\} $. Thus there is a net $\left( m_{\iota }\right)
_{\iota }$ given by
\begin{equation*}
m_{\iota }=\sum_{j=1}^{n_{\iota }}w_{\iota j}\overline{\chi (g_{\iota j})}
U_{g_{\iota j}}
\end{equation*}
with $w_{\iota 1}+...+w_{\iota n_{\iota }}=1$ and $w_{\iota j}\geq 0$, which
converges strongly to $P_{\chi }$.

Now for any $a\in M$ the net $\left( c_{\iota }\right) _{\iota }$ given by
\begin{equation*}
c_{\iota }(a)=\sum_{j=1}^{n_{\iota }}w_{\iota j}\overline{\chi (g_{\iota j})}
\theta _{g_{\iota j}}(a)
\end{equation*}
has a weakly convergent subnet, say $\left( c_{\iota (\kappa )}(a)\right)
_{\kappa }$, since the unit ball of $B(H)$ is weakly compact
\cite[Proposition 2.4.2]{BR}. Since $M$ is a von Neumann algebra, it is
weakly closed, so the limit of this subnet, say $a_{0}$, is in $M$. Since
$\left( m_{\iota (\kappa )}\right) _{\kappa }$ converges strongly to $P_{\chi
}$ and $U_{g}\Omega =\Omega $, we have that
\begin{equation*}
\left\langle x,P_{\chi }a\Omega \right\rangle =\lim_{\kappa }\left\langle
x,m_{\iota (\kappa )}a\Omega \right\rangle =\lim_{\kappa }\left\langle
x,c_{\iota (\kappa )}(a)\Omega \right\rangle =\left\langle x,a_{0}\Omega
\right\rangle
\end{equation*}
for any $x\in H$, so $P_{\chi }a\Omega =a_{0}\Omega $. It follows that
$\theta _{g}\left( a_{0}\right) \Omega =\chi (g)a_{0}\Omega $, but $\Omega $
is separating for $M$, hence $a_{0}\in M_{\chi }$.

We remark that one can in fact go further: Taking $b\in M^{\prime }$ we have
\begin{equation*}
c_{\iota (\kappa )}(a)b\Omega =bm_{\iota (\kappa )}a\Omega \rightarrow
bP_{\chi }a\Omega =a_{0}b\Omega
\end{equation*}
since $c_{\iota (\kappa )}(a)\in M$. Since $M^{\prime }\Omega $ is dense in
$H$, it is now straightforward to show that $b\Omega $ can be replaced by any
$x\in H$, in other words $\left( c_{\iota (\kappa )}(a)\right) $ actually
converges strongly to $a_{0}$. We will only need the convergence $c_{\iota
(\kappa )}(a)\Omega \rightarrow a_{0}\Omega $, though.

Now consider any $x\in H_{\chi }$ and $\varepsilon >0$, and let $a\in M$ be
such that $\left\| a\Omega -x\right\| <\varepsilon $. Then
\begin{equation*}
\left\| \overline{\chi (g)}\theta _{g}(a)\Omega -x\right\| =\left\|
U_{g}a\Omega -\chi (g)x\right\| <\varepsilon
\end{equation*}
so $\left\| c_{\iota (\kappa )}(a)\Omega -x\right\| <\varepsilon $. But by
the convergence above, there is a $\kappa $ such that $\left\| a_{0}\Omega
-c_{\iota (\kappa )}(a)\Omega \right\| <\varepsilon $, so $\left\|
a_{0}\Omega -x\right\| <2\varepsilon $. This proves that $M_{\chi }\Omega $
is dense in $H_{\chi }$ as required, since clearly $M_{\chi }\Omega \subset
H_{\chi }$.

(ii) Now we prove a number of properties of $T$. First, since
$U_{g}H_{0}\subset H_{0}$ from the definition of $H_{0}$ while $U_{g}^{\ast
}\Omega =\Omega $, we clearly have $\tau _{g}(T)=\theta_{g}(T)=T$. It
follows that $\tau _{g}\in B(T)$.

Next, set
\begin{equation*}
S:=\text{span}\bigcup_{\chi \in \mathbb{\sigma }_{\mathbf{A}}}M_{\chi }
\end{equation*}
where ``span'' means finite linear combinations. From (i) it follows that
$\overline{S\Omega }=H_{0}$. It is also readily verified that $S$ is a $\ast $
-algebra contained in $T$; in fact, $S$ is the $\ast $-algebra generated by
the eigenoperators of $\theta $. In particular $\overline{T\Omega }=H_{0}$.

Lastly we show that $T$ is weakly closed, so consider any $a$ in the weak
closure of $T$, and a net $\left( a_{\iota }\right) $ in $T$ converging
weakly to $a$. Since $M$ is weakly closed, $a\in M$. Furthermore, from the
definition of $T$, for every $x\in H_{0}^{\perp }$ we have $0=\left\langle
x,a_{\iota }\Omega \right\rangle \rightarrow \left\langle x,a\Omega
\right\rangle $, so $a\Omega \in H_{0}$. Therefore, again from the
definition of $T$, we have $a\in T$ as required.

(iii) Everything so far in the proof holds for non-abelian $G$ as well, but
in the rest of the proof we do make use of the fact that $G$ is abelian,
since we are going to work with its Bohr compactification.

In particular we now show that $T$ is a von Neumann algebra using the group
$\Gamma $ given by Lemma 5.3 applied to $T$ with $W=U$. Note that all the
requirements in Lemma 5.3 are satisfied: $M^{\prime }\subset T^{\prime }$ so
$T^{\prime }\Omega $ is dense in $H$, and for every $a\in T$ the closure
$K_{a}$ of the corresponding orbit in $H$ is indeed compact, since each
element of $H_{0}$ has a totally bounded (i.e. relatively compact) orbit
under $U$ (in fact $H_{0}$ contains all such elements \cite[Lemma 6.6]{BDS},
but this fact will only be used in (iv) below).

Assign the discrete topology to $G$. Since $G$ is abelian, we can imbed it
into its Bohr compactification $\bar{G}$. Let us denote this canonical
imbedding by $\iota :G\rightarrow \bar{G}$. It is an injective group
homomorphism with $\iota (G)$ dense in $\bar{G}$. The dual map of this is a
group isomorphism $\hat{\iota}:\widehat{\bar{G}}\rightarrow \widehat{G}$.
For $\chi\in\widehat{G}$, set
$\left\langle\cdot,\chi\right\rangle:=\hat{\iota}^{-1}(\chi)$,
then it follows that
$\left\langle \iota (g),\chi\right\rangle=\chi(g)$.
Note that since $\tau _{g}\in B(T)$ is an isometry,
it follows that $\tau :G\rightarrow \Gamma :g\mapsto \tau _{g}$ is a
well-defined group homomorphism (and it is continuous, since $G$ is
discrete) which by the universal property of the Bohr compactification
\cite[Proposition (4.78)]{F} can be extended to a continuous group
homomorphism
\begin{equation*}
\bar{\tau}:\bar{G}\rightarrow \Gamma :g\mapsto \bar{\tau}_{g}
\end{equation*}
which means that $\bar{\tau}_{\iota (g)}=\tau _{g}$ for all $g\in G$.

Since $\Gamma \rightarrow H:\Theta \mapsto \Theta (a)y$ is continuous, so is
$\bar{G}\rightarrow H:g\mapsto \bar{\tau}_{g}(a)y$ for all $a\in T$ and
$y\in H$. Denoting the normalized Haar measure on the compact group $\bar{G}$
by $m$, we can for every $a\in T$ and $\chi \in \widehat{G}$ define a unique
$a_{\chi}\in B(H)$ by requiring
\begin{equation*}
\left\langle x,a_{\chi }y\right\rangle =\int_{\bar{G}}\overline{\left\langle
g,\chi \right\rangle }\left\langle x,\bar{\tau}_{g}(a)y\right\rangle dm(g)
\end{equation*}
for all $x,y\in H$. If $f$ is a weakly continuous linear functional on $B(H)$
it follows (see for example \cite[Theorem V.3.9]{DS}) that it is a finite
linear combination of such $\left\langle x,(\cdot )y\right\rangle $ forms,
hence
\begin{equation}
f(a_{\chi })=\int_{\bar{G}}\overline{\left\langle g,\chi \right\rangle }
f\left( \bar{\tau}_{g}(a)\right) dm(g)  \tag{1}
\end{equation}
so if $f$ is zero on $T$ we have $f(a_{\chi })=0$ and therefore by
\cite[Corollary V.3.12]{DS} $a_{\chi }\in T$ since $T$ is weakly closed. For
$g\in G$ it follows from the definition of $a_{\chi}$ that
\begin{eqnarray*}
\left\langle x,\tau _{g}(a_{\chi })y\right\rangle
&=&\left\langle U_{g}^{\ast }x,a_{\chi }U_{g}^{\ast }y\right\rangle
=\int_{\bar{G}}\overline{\left\langle h,\chi \right\rangle }\left\langle U_{g}^{\ast }x,\bar{\tau}
_{h}(a)U_{g}^{\ast }y\right\rangle dm(h) \\
&=&\int_{\bar{G}}\overline{\left\langle h,\chi \right\rangle }\left\langle
x, \bar{\tau}_{\iota (g)h}(a)y\right\rangle dm(h) \\
&=&\int_{\bar{G}}\overline{\left\langle \iota (g)^{-1}h,\chi \right\rangle }
\left\langle x,\bar{\tau}_{h}(a)y\right\rangle dm(h) \\
&=&\overline{\left\langle \iota (g)^{-1},\chi \right\rangle }\int_{\bar{G}}
\overline{\left\langle h,\chi \right\rangle }\left\langle x,\bar{\tau}
_{h}(a)y\right\rangle dm(h) \\
&=&\chi (g)\left\langle x,a_{\chi }y\right\rangle
\end{eqnarray*}
and therefore $\tau _{g}(a_{\chi })=\chi (g)a_{\chi }$ which means that
$a_{\chi }\in S$, since $a_{\chi }\in T\subset M$. We now use this result to
show that $S$ is weakly dense in $T$, so take any $a\in T$. Let $f$ be any
weakly continuous linear functional on $B(H)$ which is zero on the weak
closure $\overline{a_{\widehat{G}}}$ of the span of $\left\{ a_{\chi }:\chi
\in \widehat{G}\right\} $. From (1) and the Plancherel theorem it follows
that $g\mapsto f\left( \bar{\tau}_{g}(a)\right) $ is zero in $L^{2}(\bar{G})$.
But this function is continuous by arguments as above, so
$f\left(\bar{\tau }_{g}(a)\right) =0$ for all $g\in\bar{G}$, for if not,
then $f$ would be non-zero on an open neighbourhood of some $g\in \bar{G}$,
contradicting the fact that the Haar measure of a non-empty open set is non-zero.
In particular for $g=1$ we find $f(a)=0$. This means that
$a\in \overline{a_{\widehat{G}}}$ and therefore $S$ is weakly dense in $T$.
Since $S$ is a $\ast$-algebra containing the identity operator, we conclude that $T$ is a
von Neumann algebra by the von Neumann density theorem
\cite[Corollary 2.4.15]{BR}.

(iv) By (ii) and (iii) we know that $\mathbf{T}$ is indeed a subsystem of
$\mathbf{A}$. It has $\sigma _{\mathbf{T}}=\sigma _{\mathbf{A}}$, since it
can be cyclically represented on $H_{0}$ (because of $\overline{T\Omega
_{\mu }}=H_{0}$) with the dynamics given by $U_{g}|_{H_{0}}$. Since $S$
above is the $\ast $-algebra generated by the eigenoperators of $\theta $,
we have $S\subset \pi _{\mu }(A^{K})\subset T$. Hence $T=\pi _{\mu }(A^{K})$,
since $\pi _{\mu }(A^{K})$ is a von Neumann algebra and therefore weakly
closed. From this the isomorphism between $\mathbf{T}$ and $\mathbf{A}^{K}$
follows. If $\mathbf{F}$ is any compact subsystem of $\mathbf{A}$, and $a\in
F$, then by definition of compactness of a system the orbit of $\gamma _{\mu
}(\zeta (a))$ under $U$ is totally bounded (keep in mind that $\left(
F,\lambda \right) $ can be cyclically represented on $\overline{\gamma _{\mu
}(\zeta (F))}$, with the corresponding unitary representation of $\varphi $
given by the restriction of $U_{g}$ to this space), thus $\pi _{\mu }(\zeta
(a))\Omega _{\mu }\in H_{0}$ according to \cite[Lemma 6.6]{BDS}. In other
words $\pi _{\mu }(\zeta (a))\in T$, and therefore $\zeta (a)\in A^{K}$
which means that $\mathbf{A}^{K}$ is the largest compact subsystem of
$\mathbf{A}$.
\end{proof}

Note that the sole reason we needed to assume that $G$ is abelian in this
theorem, is that we use its Bohr compactification in the proof.

Whereas $\mathbf{A}^{\alpha }$ is by definition the largest identity
subsystem of $\mathbf{A}$, we have seen above that $\mathbf{A}^{K}$ is the
largest compact subsystem of $\mathbf{A}$, at least when $G$ is abelian.
Clearly $A^{\alpha }\subset A^{K}$. We mention that a very simple version of
this in the context of noncommutative topological dynamics was also
discussed in \cite[Section 2]{LL} and \cite[Definition 1.2]{A}. In analogy
to Proposition 4.2, Theorem 5.4 says that
$H_{0}=\overline{\gamma _{\mu}(A^{K})}$.

Before we return to relatively independent joinings we note the following
generalization of \cite[Theorem 6.8]{BDS} (also see \cite[Proposition 5.4]{NSZ}).
Note that in terms of Definition 5.1's notation, we call $\mathbf{A}$
\textit{weakly mixing} if $\dim H_{0}=1$.

\begin{corollary}
A system $\mathbf{A}$ with $G$ assumed abelian, is weakly mixing if and only
if $\mathbf{A}^{K}$ is the trivial system.
\end{corollary}

\begin{proof}
If $\mathbf{A}^{K}$ is trivial, then so is $\mathbf{T}$ in Theorem 5.4, so
$H_{0}=\overline{T\Omega _{\mu }}=\mathbb{C}\Omega _{\mu }$. Conversely, if
$\mathbf{A}$ is weakly mixing, then $\overline{T\Omega _{\mu }}=\mathbb{C}
\Omega _{\mu }$, but $\Omega _{\mu }$ is separating for $\pi _{\mu }(A)$ and
therefore for $T$, so $T=\mathbb{C}$, hence $\mathbf{A}^{K}$ is trivial.
\end{proof}

In addition to the assumptions in this corollary, in \cite[Theorem 6.8]{BDS}
it was also assumed that $\mathbf{A}$ is ergodic.

Now we finally return to relatively independent joinings, namely an analogue
of Theorem 4.3. The proof is very similar to that of Theorem 4.3.

\begin{theorem}
Assume that $G$ is abelian, and let $\mathbf{F}$ be a compact modular
subsystem of $\mathbf{A}$. If $\mathbf{A}$ and $\mathbf{\tilde{B}}$ are
disjoint over $\mathbf{F}$ for $\mathbf{B}=\mathbf{A}^{K}$, then $\mathbf{F}$
is isomorphic to $\mathbf{A}^{K}$. On the other hand, if $\mathbf{F}$ is
isomorphic to $\mathbf{A}^{K}$, then $\mathbf{A}$ and $\mathbf{\tilde{B}}$
are disjoint over $\mathbf{F}$ for every compact system $\mathbf{B}$ which
has $\mathbf{F}$ as a modular subsystem.
\end{theorem}

\begin{proof}
Assume without loss that $\mathbf{F}$ is imbedded in $\mathbf{A}$ by
inclusion. Suppose $\mathbf{F}$ is not isomorphic to $\mathbf{A}^{K}$, then
$F$ is strictly contained in $A^{K}$ by Theorem 5.4. Note that $\mathbf{F}$
is a modular subsystem of $\mathbf{A}^{K}$, since $\mathbf{F}$ and
$\mathbf{A}^{K}$ are modular subsystems of $\mathbf{A}$. Now apply Proposition 3.5
with $\mathbf{B}=\mathbf{R}=\mathbf{A}^{K}$ to obtain $\mu \odot _{\psi }
\tilde{\nu}\neq \mu \odot _{\lambda }\tilde{\nu}$. But both these joinings
are contained in $J_{\lambda }\left( \mathbf{A},\mathbf{\tilde{B}}\right) $
by Propositions 3.2 and 3.5, so $\mathbf{A}$ and $\mathbf{\tilde{B}}$ are
not disjoint over $\mathbf{F}$.

Conversely, assume that $\mathbf{F}$ is isomorphic (and therefore equal) to
$\mathbf{A}^{K}$. Consider any
$\omega \in J_{\lambda }\left( \mathbf{A},\mathbf{\tilde{B}}\right) $.
For any eigenvector $y\in H_{\tilde{\nu}}$ of $\tilde{V}$ with eigenvalue
$\chi \in \sigma _{\mathbf{\tilde{B}}}$ we have
\begin{equation*}
U_{g}P_{\omega }y=P_{\omega }\tilde{V}_{g}y=\chi (g)P_{\omega }y
\end{equation*}
with $P_{\omega }$ as in the proof of Theorem 4.3. So $P_{\omega }y=0$ or
$\chi \in \sigma _{\mathbf{A}}$. Therefore $P_{\omega }y\in H_{0}=H_{\lambda
} $ with $H_{0}$ as in Definition 5.1 and $H_{\lambda }$ as given in the
setting described before Proposition 3.6. Since $\mathbf{\tilde{B}}$ has
discrete spectrum (i.e. $H_{\tilde{\nu}}$ is spanned by the eigenvectors of
$\tilde{V}$) as remarked in \cite[Proposition 2.6]{D2} (also see
\cite[Section 2.4]{K} for the general theory), it follows that $P_{\omega
}y\in H_{\lambda }$ for all $y\in H_{\tilde{\nu}}$. Since $P_{\omega }$ is
the projection of $H_{\omega }$ onto $H_{\mu }$ restricted to $H_{\tilde{\nu}}$,
it follows that $H_{\tilde{\nu}}\perp \left( H_{\mu }\ominus H_{\lambda}\right)$,
and hence $\omega =\mu \odot _{\lambda }\tilde{\nu}$ by Proposition 3.6.
\end{proof}

In the case of trivial $\mathbf{F}$ this theorem gives a variation on
\cite[Theorems 2.7 and 2.8]{D2} in that the latter assumed ergodicity of the
systems involved, but did not require $G$ to be abelian, namely:

\begin{corollary}
Assuming $G$ is abelian, $\mathbf{A}$ is weakly mixing if and only if it is
disjoint from all compact systems.
\end{corollary}

\begin{proof}
Use trivial $\mathbf{F}$ in Theorem 5.6, together with Proposition 3.4 and
Corollary 5.5.
\end{proof}

As in Section 4, these results illustrate how relatively independent
joinings relate to the structure of W*-dynamical systems, in particular to
certain types of subsystems.

\section{Ergodicity and compactness}

In this concluding section we turn away from joinings and focus on
subsystems. The goal is to study $\mathbf{A}^{K}$, which was defined in
Proposition 5.2, when $\mathbf{A}$ is ergodic, i.e. when $\mathbf{A}^{\alpha}$
is trivial. This is related to, and largely motivated by \cite[Theorem 1.3]{JP}
(see Corollary 6.2 below).

Remember that a state $\mu $ on a von Neumann algebra $A$ is called
\emph{tracial} if $\mu (ab)=\mu (ba)$ for all $a,b\in A$. The following result
strengthens the fact that $\mathbf{A}^{K}$ is a modular subsystem of
$\mathbf{A}$ when the latter is ergodic, and note that in this result we need
not assume $G$ is abelian:

\begin{theorem}
Let $\mathbf{A}$ be ergodic. Then $A^{K}$ is contained in the fixed point
algebra of the modular group $\sigma ^{\mu }$ of $\mu $. Furthermore, $\mu
^{K}$ is then tracial, and $A^{K}$ therefore a finite von Neumann algebra.
\end{theorem}

\begin{proof}
Let $u$ be any eigenoperator of $\alpha $, with eigenvalue $\chi $. Then as
in the proof of Proposition 4.2
\begin{equation*}
\alpha _{g}(\sigma _{t}^{\mu }(u))=\sigma _{t}^{\mu }(\alpha _{g}(u))=\chi
(g)\sigma _{t}^{\mu }(u).
\end{equation*}
Since $\mathbf{A}$ is ergodic and $u$ and $\sigma _{t}^{\mu }(u)$ are
eigenoperators with the same norm and the same eigenvalue, it follows from
\cite[Lemma 2.1(3)]{S} that there is a number, say $\Lambda (t)\in \mathbb{C}$,
with $\left| \Lambda (t)\right| =1$, such that
\begin{equation*}
\sigma _{t}^{\mu }(u)=\Lambda (t)u
\end{equation*}
for all $t\in \mathbb{R}$. Now we show that $\sigma _{t}^{\mu }(u)=u$ (also
see the proof of \cite[Lemma 5.2]{D2}). It follows from the group property
of $\sigma ^{\mu }$ that $\Lambda (s+t)=\Lambda (s)\Lambda (t)$ for all
$s,t\in \mathbb{R}$, and since $t\mapsto \left\langle x,\pi _{\mu }(\sigma
_{t}^{\mu }(u))y\right\rangle $ is continuous for all $x,y\in H_{\mu }$, it
follows that $t\mapsto \Lambda (t)$ is continuous. Therefore
\begin{equation*}
\Lambda (t)=e^{i\theta t}
\end{equation*}
for all $t\in \mathbb{R}$ for some $\theta \in \mathbb{R}$; see for example
\cite[p. 12]{R}. Denoting the modular operator associated with $\left( \pi
_{\mu }(A),\Omega _{\mu }\right) $ by $\Delta $, it follows that $\Delta
^{it}\pi _{\mu }(u)\Omega _{\mu }=\pi _{\mu }(\sigma _{t}^{\mu }(u))\Omega
_{\mu }=e^{i\theta t}\pi _{\mu }(u)\Omega _{\mu }$, hence by the definition
of $J_{\mu }\Delta ^{1/2}$ (see for example \cite[Section 2.5.2]{BR})
\begin{equation*}
J_{\mu }\pi _{\mu }(u)^{\ast }\Omega _{\mu }=J_{\mu }\left( J_{\mu }\Delta
^{1/2}\right) \pi _{\mu }(u)\Omega _{\mu }=\Delta ^{1/2}\pi _{\mu }(u)\Omega
_{\mu }=e^{\theta /2}\pi _{\mu }(u)\Omega _{\mu }
\end{equation*}
and by taking the norm both sides we conclude that $e^{\theta /2}=1$, since
$\alpha _{g}(uu^{\ast })=\left| \chi (g)\right| ^{2}uu^{\ast }=uu^{\ast }$
and $\alpha _{g}(u^{\ast }u)=u^{\ast }u$, implying $uu^{\ast }=u^{\ast }u\in
\mathbb{C}1\backslash \{0\}$ by ergodicity and the fact that $uu^{\ast }$
and $u^{\ast }u$ are both positive. Therefore $\theta =0$. This proves that
\begin{equation*}
\sigma _{t}^{\mu }(u)=u
\end{equation*}
for all $t\in \mathbb{R}$. Then by the definition of $A^{K}$ in Proposition
5.2, it is indeed contained in the fixed point algebra of the modular group
$\sigma ^{\mu }$, since the latter algebra is a von Neumann algebra
containing all eigenoperators of $\alpha $, as just shown.

Since $\mathbf{A}^{K}$ is a modular subsystem of $\mathbf{A}$, the modular
group $\sigma ^{\mu ^{K}}$ of $\mu ^{K}$ is simply given by the restriction
of $\sigma _{t}^{\mu }$ to $A^{K}$. Hence $\sigma _{t}^{\mu ^{K}}=$ id$_{A^{K}}$
for all $t$ by what we have shown above. It follows that the
modular operator associated with the cyclic representation of $\left(
A^{K},\mu ^{K}\right) $ is the identity operator. Therefore by the basic
definitions of Tomita-Takesaki theory, and writing the cyclic representation
of $\left( A^{K},\mu ^{K}\right) $ as $\left( H,\pi ,\Omega \right) $ and
the corresponding modular conjugation as $J$, we have $Ja\Omega =a^{\ast
}\Omega $ for all $a\in \pi (A^{K})$. From this it follows that $\mu ^{K}$
is tracial, namely for any $a,b\in A^{K}$ we have
\begin{eqnarray*}
\mu (ab) &=&\left\langle \Omega ,\pi (ab)\Omega \right\rangle =\left\langle
\pi (a)^{\ast }\Omega ,\pi (b)\Omega \right\rangle =\left\langle J\pi
(a)\Omega ,\pi (b)\Omega \right\rangle \\
&=&\left\langle J\pi (b)\Omega ,\pi (a)\Omega \right\rangle =\left\langle
\Omega ,\pi (ba)\Omega \right\rangle \\
&=&\mu (ba)
\end{eqnarray*}
and the existence of a faithful (numerical) trace on $A^{K}$ implies that it
is a finite von Neumann algebra (see for example \cite[p. 505]{KR2}).
\end{proof}

To clarify the connection with \cite[Theorem 1.3]{JP}, we deduce it as a
corollary from this theorem (but we only state it for an inverse temperature
of $-1$):

\begin{corollary}
Suppose that the dynamics $\alpha $ of $\mathbf{A}$ is given by the modular
group $\sigma ^{\mu }$, i.e. $G=\mathbb{R}$ and $\alpha _{t}=\sigma
_{t}^{\mu }$ for all $t\in \mathbb{R}$. If $\mathbf{A}$ is ergodic, it then
follows that it is weakly mixing.
\end{corollary}

\begin{proof}
By Theorem 6.1 $A^{K}\subset A^{\alpha }$, but $A^{\alpha }=\mathbb{C}$,
since $\mathbf{A}$ is ergodic, therefore $A^{K}=\mathbb{C}$. Since $\mathbb{R}$
is an abelian group, it follows from Corollary 5.5 that $\mathbf{A}$ is
weakly mixing.
\end{proof}

One may wonder whether more generally, when $\mathbf{A}^{\alpha }$ is not
necessarily trivial, one has $A^{K}=A^{\alpha }$ if $\alpha $ is the modular
group of $\mu $. However on $A=B(H)$ with $H$ a finite dimensional Hilbert
space, one has the simple counter example given by the Gibbs state (at
inverse temperature $-1$) $\mu (a)=$ Tr$(e^{h}a)/$Tr$(e^{h})$ where $h\in A$
is hermitian, in which case $\alpha _{t}(a)=e^{iht}ae^{-iht}$ is the modular
group of $\mu $; the system $\mathbf{A}$ so obtained is compact but not in
general an identity system, so $A^{K}=A\neq A^{\alpha }$.

Another corollary of Theorem 6.1 is the following (also see
\cite[Proposition 7.2]{OPT}):

\begin{corollary}
If a system $\mathbf{A}$ is ergodic and compact, and $G$ is abelian, then
$\mu $ is necessarily tracial and $A$ a finite von Neumann algebra.
\end{corollary}

\begin{proof}
By Theorem 5.4 $\mathbf{A}^{K}=\mathbf{A}$, and the result then follows from
Theorem 6.1.
\end{proof}

Results of the nature of Corollary 6.3 have of course been well studied. However, our
assumption that the system is compact makes Corollary 6.3 a relatively easy
result. In the literature much more difficult results have been obtained;
see \cite{H, S67, S71, S, L, AH, OPT, HLS} for the development.

The results of this section illustrate how subsystems can be applied to
derive interesting properties of systems, and show how aspects of the work
above fit into previous literature on noncommutative dynamical systems.

\section*{Acknowledgments}

I thank Anton Str\"{o}h for useful discussions, and Richard de Beer for
pointing the paper \cite{JP} out to me long ago. This work was supported by
the National Research Foundation of South Africa.

\end{document}